\documentclass[11pt]{amsart}
\usepackage{amssymb}
\usepackage{amsmath}
\usepackage[latin1]{inputenc}

\usepackage[english]{babel}

\newtheorem{theorem}{Theorem}[section]
\newtheorem{lemma}[theorem]{Lemma}

\newtheorem{corollary}[theorem]{Corollary}
\newtheorem{e-definition}[theorem]{Definition\rm}

\begin{document}


\selectlanguage{english}
\title[Essential norm of composition operators on Bloch-type spaces]{An estimation for the essential norm of composition operators acting on Bloch-type spaces}


\selectlanguage{english}

\author[J. C. Ramos-Fern\'andez]{Julio C. Ramos Fern\'andez}

\address[Julio C. Ramos-Fern\'andez]{Departamento de Matem\'atica, Universidad de Oriente, 6101 Cuman\'a,
Edo. Sucre, Venezuela}

\begin{abstract}
Let  $\mu$ be any weight  function defined on the unit disk $\Bbb D$ and let $\phi$ be an analytic self-map of  $\Bbb D$.
 In the present paper we show that the essential norm of composition operator $C_\phi$ mapping from the $\alpha$-Bloch space, with $\alpha >0$,  to $\mu$-Bloch space $\mathcal{B}^\mu$ is comparable to
 $$
  \limsup_{|a|\to 1^-}\left\|\sigma_a\circ\phi\right\|_{\mathcal{B}^{\mu}},
 $$
 where, for $a\in\Bbb D$, $\sigma_a$ is a certain special function in $\alpha$-Bloch space.

\vspace{0.1cm}


\end{abstract}
\maketitle

\selectlanguage{english}

\section{Introduction}
Let $\mu$ be denotes what we call a {\em weight} on the unit disk ${\Bbb D}$ of the complex plane $\Bbb C$; that is, $\mu$ is a
bounded, continuous and strictly positive function defined on ${\Bbb
D}$, and let $H(\Bbb D)$ be the space of all holomorphic functions on $\Bbb D$, which is equipped with the topology of uniform convergence on compact subsets of $\Bbb D$. The {\em $\mu$-Bloch space} $\, \mathcal{B}^\mu({\Bbb D})$,
which we denote more briefly by ${\mathcal B}^{\mu}$, consists of all
$f\in H({\Bbb D})$ such that
$$
 \left\|f\right\|_{\mu}:=\sup_{z\in{\Bbb
 D}}\mu\left(z\right)\left|f'(z)\right|<\infty.
$$
$\mu$-Bloch spaces are called {\em weighted Bloch spaces}.  For weights $\mu$ on ${\Bbb D}$, a Banach space
structure  on ${\mathcal B}^{\mu}$ arises if it is given the norm
$$\|f\|_{\mathcal{B}^\mu}:=\left|f(0)\right|+||f||_\mu.$$
These Banach spaces provide a natural setting in which one can study
properties of various operators. For instance,  Attele in \cite{At92}
proved that if $\mu_1(z):=w(z)\log\frac{2}{w(z)}$, where
$w(z):=1-|z|^2$ and $z\in\Bbb D$, then the Hankel operator $H_f$
induced by a function $f$ in the Bergman space $A^2({\Bbb D})$ (see
\cite[~Ch.~2]{CowMac95}) is bounded if and only if $\,f\in
B^{\mu_1}$, thus giving one reason, and not the only
reason, why log-Bloch-type spaces are of interest. When $\mu(z)=v_\alpha (z) :=\left(1-|z|^2\right)^\alpha$ with $\alpha>0$ fixed, then we get back the $\alpha$-Bloch space which is denoted as $\mathcal{B}^\alpha$ and when $\alpha=1$ we obtain the Bloch space $\mathcal{B}$.

\vspace{0.3cm}
A holomorphic function $\phi$ from the unit disk $\Bbb D$ into itself induces a linear operator $C_\phi$, defined by $C_\phi(f)=f\circ \phi$, where $f\in H(\Bbb D)$.  $C_\phi$ is called the \textit{composition operator} with
\textit{symbol} $ \, \phi. \, $  Composition operators continue to
be widely studied on many subspaces of $\, H({\Bbb D})\, $ and
particularly in Bloch-type spaces.

\vspace{0.3cm}
The study of the properties of composition operators on Bloch-type spaces began with the celebrated work of Madigan and Matheson in \cite{MM95}, where they characterized the continuity and compactness of composition operators acting on the Bloch space $\mathcal{B}$. Many extensions of the Madigan and Matheson's results have appeared (see for instance \cite{Ra11EM} and  a lot of references therein). In particular, Xiao in \cite{Xi01} has extended the results by Madigan and Matheson in \cite{MM95} to composition operators $C_\phi$ acting between $\alpha$-Bloch spaces. Recently, many authors have found new criteria for the continuity and compactness of composition operators acting on Bloch-type spaces en terms of the $n$-th power of the symbol  $\phi$ and the norm of the  $n$-th power of the identity function on  $\Bbb D$. The first result of this kind appears in  2009 and it is due to   Wulan, Zheng, and Zhu (\cite{WZheZhu09}), in turn, their result was extended to $\alpha$-Bloch spaces by Zhao in \cite{Zhao10}. Another criterion for the continuity and compactness of composition operators on Bloch space is due to Tjani  in \cite{Tj03} (see also \cite{Wu07} or more recently  \cite{WZheZhu09}), she showed the following result:

\begin{theorem}[\cite{Tj03}]\label{th-Tjani}
 The composition operator $C_\phi$ is compact on $\mathcal{B}$  if and only if \, $\phi\in\mathcal{B}$ and
   $$
     \lim_{\left|a\right|\to 1^-}\left\|\varphi_a\circ \phi\right\|_{\mathcal{B}}=0,
   $$
   where  $\varphi_a$ is a M\"{o}bius transformation from the unit disk onto itself;  that is,
$ \varphi_a(z)=\left(a-z\right)/\left(1-\overline{a}z\right)$, with $z\in\Bbb D$.
\end{theorem}
\noindent This last result has been recently extended to $\alpha$-Bloch spaces by Malav\'e and Ramos-Fern\'andez in \cite{MR13}.

\vspace{0.2cm}
The essential norm of a continuous linear operator
$T$ between normed linear spaces $X$ and $Y$ is its distance from
the compact operators; that is, $\|T\|^{X\rightarrow Y}_e =
\inf\left\{\|T - K\|^{X\rightarrow Y} : K: X\to Y\, \text{ is compact}\right\}$, where
$\|\cdot\|^{X\rightarrow Y}$ denotes the operator norm. Notice that
$\|T\|^{X\rightarrow Y}_e = 0$ if and only if $T$ is compact, so
that estimates on $\|T\|_e^{X\rightarrow Y}$ lead to conditions for
$T$ to be compact.
The essential norm of a composition operator on ${\mathcal B}$ was calculated by A. Montes-Rodr\'{\i}guez in \cite{MoRo99}.
He obtained similar results for essential norms of weighted composition
operators between weighted Banach spaces of analytic functions in \cite{MoRo00}.  Other results in this direction appear in the paper by  Contreras and  Hern\'andez-D\'{\i}az in \cite{ContH00}; in
particular, formulas for the essential norm of  weighted composition
operators on the $\alpha$-Bloch spaces of ${\Bbb B}_n$ were obtained
(see also the paper of MacCluer and Zhao \cite{MacZhao03}).
Recently, many extensions of the above results have appeared in the
literature; for instance, the reader is referred to the paper of
 Yang and  Zhou \cite{YaZho10} and several references therein. Zhao in
\cite{Zhao10} gave a formula for the essential norm of
$C_\phi:\mathcal{B}^\alpha\to\mathcal{B}^\beta$
in terms of an expression involving norms of powers of $\phi$. More
precisely, he showed that
 $$
  \left\|C_\phi\right\|_e^
  {\mathcal{B}^\alpha\to\mathcal{B}^\beta}=\left( \frac{e}{2\alpha} \right) ^\alpha\limsup_{j\to\infty}j^{\alpha-1}\left\|\phi^j\right\|_{{\mathcal B}^{\beta}}.
 $$
It follows from the discussion at the beginning of this paragraph
that $C_\phi:\mathcal{B}^\alpha\to\mathcal{B}^\beta$ is compact if and only if
 $$
  \lim_{j\to\infty}j^{\alpha-1}\left\|\phi^j\right\|_{{\mathcal B}^{\beta}}=0.
 $$
The Zhao's results in \cite{Zhao10} have been extended recently to the weighted Bloch spaces by Castillo, Clahane, Far\'{\i}as and Ramos-Fern\'andez in \cite{CCFR13}. Also, Hyv\"{a}rinen,  Kemppainen,  Lindstr\"{o}m,  Rautio and Saukko in \cite{gblochpaper} obtained necessary and
sufficient conditions for boundedness and an expression
characterizing the essential norm of a weighted composition operator
between general weighted Bloch spaces ${\mathcal B}^{\mu}$, under
the technical requirements that $\mu$ is radial, and that it is
non-increasing and tends to zero toward the boundary of ${\Bbb D}$.

\vspace{0.2cm}
The goal of the present paper is to give an estimate of the essential norm of composition $C_\phi$ mapping from $\mathcal{B}^\alpha$ to $\mathcal{B}^\mu$ which implies Theorem \ref{th-Tjani} and the result given by Malav\'e and Ramos-Fern\'andez in \cite{MR13}. More precisely, in the next section we will show the following result.
\begin{theorem}\label{th-principal}
 Let $\phi$ be an analytic self-map of the unit disk $\Bbb D$. Then for the essential norm of the composition operator $C_\phi:\mathcal{B}^\alpha\to\mathcal{B}^{\mu}$ we have
 \begin{equation}\label{ed-equivalentes}
  \left\|C_\phi\right\|_e^{\mathcal{B}^\alpha\to\mathcal{B}^\mu}\sim\limsup_{|a|\to 1^-}\left\|\sigma_a\circ\phi\right\|_{\mathcal{B}^{\mu}}.
 \end{equation}
\end{theorem}
\noindent  The relation (\ref{ed-equivalentes}) means that there is
a positive constant $M_\alpha$, depending only on $\alpha$, such that
$$
 \frac{1}{M_\alpha}\left\|C_\phi\right\|_e^{\mathcal{B}^\alpha\to\mathcal{B}^\mu}\leq \limsup_{|a|\to 1^-}\left\|\sigma_a\circ\phi\right\|_{\mathcal{B}^{\mu}}\leq M_\alpha\, \left\|C_\phi\right\|_e^{\mathcal{B}^\alpha\to\mathcal{B}^\mu}
 $$
 and the functions $\sigma_a$ with $a\in\Bbb D$ will be defined at the begin of the next section.

\section{Proof of Theorem \ref{th-principal}}
The key to our results lies in considering the following family of functions. For  $a \in \mathbb{D}$ fixed, we define
$$
\sigma_{a}(z)=(1-|a|)\left((1-\overline{a}z)^{-\alpha}-1\right),~~~~(z \in \mathbb{D}).
$$
Clearly, for each $a\in\Bbb D$, the function $\sigma_a$ has bounded derivative and for this reason we have that $\sigma_a\in \mathcal{B}^\alpha$. In fact, it is easy to see that
$$
 \sup_{a\in\Bbb D}\left\|\sigma_a\right\|_{\mathcal{B}^\alpha}\leq \alpha\,2^\alpha.
$$
Furthermore, it is clear that if $\frac{1}{2}<|a|<1$, then
\begin{equation}\label{des-cota-inf-der-sigma}
 \left|\sigma_{a}'(a)\right|\geq \frac{\alpha}{4\left(1-|a|^2\right)^{\alpha}}.
\end{equation}
Also, we can see that $\sigma_a$ goes to zero uniformly on compact subsets of $\Bbb D$ as $|a|\to 1^-$. Also, we will need the following lemma which is well known and is consequence of a more general result due to Tjani in \cite{Tj}:

\begin{lemma}\label{le-Tjani}
 The composition operator $C_\phi:\mathcal{B}^{\alpha}\to \mathcal{B}^{\mu}$ is compact if and only if given a bounded sequence
  $\left\{f_n\right\}$ in $\mathcal{B}^{\alpha}$ such that $\, f_n\to 0\, $ uniformly on compact
   subsets of $\, \Bbb D$, then $\left\|C_\phi(f_n)\right\|_{\mathcal{B}^\mu}\to 0$ as $n\to\infty$.
\end{lemma}

\vspace{0.2cm}
Now we can show Theorem \ref{th-principal}.

\noindent {\it Proof of Theorem \ref{th-principal}.} We set
$$
 L=\limsup_{|a|\to 1^-}\left\|\sigma_a\circ\phi\right\|_{\mathcal{B}^{\mu}}.
$$
Let $K:\mathcal{B}^\alpha\to \mathcal{B}^\mu$ be  any compact operator,  $a\in\Bbb D$ fixed and define
$$
 f_a(z)=\frac{1}{\alpha 2^\alpha}\sigma_a(z),\hspace{0.3cm}(z\in\Bbb D).
$$
Then $f_a$ goes to zero uniformly on compact subsets of $\Bbb D$ as $|a|\to 1^-$, $\left\|f_a\right\|_{\mathcal{B}^\alpha}\leq 1$ for all $a\in\Bbb D$ and
$$
 \left\|C_\phi - K\right\|^{\mathcal{B}^\alpha\to\mathcal{B}^\mu}\geq \left\|\left(C_\phi - K\right)f_a\right\|_{\mathcal{B}^\mu}\geq
 \frac{1}{\alpha 2^\alpha}\left\|\sigma_a\circ \phi\right\|_{\mathcal{B}^\mu} - \left\|K f_a\right\|_{\mathcal{B}^\mu}.
$$
Hence, taking $\limsup_{|a|\to 1^-}$ and using Lemma \ref{le-Tjani}, we obtain
\begin{equation}
  \left\|C_\phi\right\|_e^{\mathcal{B}^\alpha\to\mathcal{B}^\mu}\geq \frac{1}{\alpha 2^\alpha}\limsup_{|a|\to 1^-}\left\|\sigma_a\circ \phi\right\|_{\mathcal{B}^\mu}.
\end{equation}

Now, we go to show that there exists a constant $M_\alpha>0$, depending only on $\alpha$, such that
$$
  \left\|C_\phi\right\|_e^{\mathcal{B}^\alpha\to\mathcal{B}^\mu}\leq M_\alpha \limsup_{|a|\to 1^-}\left\|\sigma_a\circ \phi\right\|_{\mathcal{B}^\mu}.
$$
Bearing this in mind, we define, for $r\in [0,1]$,  the linear {\em dilation operator} $K_r:
H({\Bbb D})\to H({\Bbb D})$ by $K_r\, f=f_r$, where $f_r$,
for each $f\in  H({\Bbb D})$, is given by $f_r(z)=f(rz)$.  It is clear that if $f\in H(\Bbb D)$ then $rf_r\to f$ uniformly on compact subsets of $\Bbb D$ as $r\to 1^-$. Also, the following statements hold:
\begin{enumerate}
 \item For $r\in [0,1)$, the operator  $K_r$ is compact  on $\mathcal{B}^\alpha$,

 \item for each $r\in [0,1]$
  $$
    \left\|K_r\right\|^{{\mathcal B}^{\alpha}\rightarrow {\mathcal B}^{\alpha}}\leq 1.
  $$
\end{enumerate}
Hence, if we consider a sequence $\left\{r_n\right\}\subset (0,1)$ such that $r_n\to 1$ as $n\to\infty$ and define $K_n=K_{r_n}$, then for all $n\in\Bbb N$, the operator $C_\phi\,K_n$ is a compact  from $\mathcal{B}^\alpha$ into $\mathcal{B}^\mu$ and by definition of the essential norm we have
$$
 \left\|C_\phi\right\|_e^{\mathcal{B}^\alpha\to\mathcal{B}^\mu}\leq \limsup_{n\to\infty}\left\|C_\phi - C_\phi K_n\right\|^{\mathcal{B}^\alpha\to\mathcal{B}^\mu}.
$$
Thus, we have to show that
$$
 \limsup_{n\to\infty}\left\|C_\phi - C_\phi K_n\right\|^{\mathcal{B}^\alpha\to\mathcal{B}^\mu}\leq M_\alpha L.
$$
To see this, consider any $f\in \mathcal{B}^\alpha$ such that $\left\|f\right\|_{\mathcal{B}^\alpha}=1$, then since
$$
 \left\|\left(C_\phi - C_\phi K_n\right)f\right\|_{\mathcal{B}^{\mu}}=\left|f\left(\phi(0)\right)-f\left(r_n\phi(0)\right)\right|+ \left\|\left(f-f_{r_n}\right)\circ \phi\right\|_{\mu}
$$
and $\left|f\left(\phi(0)\right)-f\left(r_n\phi(0)\right)\right|\to 0$ as $n\to\infty$, it is enough to show that
$$
 \limsup_{n\to\infty}\left\|\left(f-f_{r_n}\right)\circ \phi\right\|_{\mu}\leq M_\alpha L.
$$
Furthermore, since $r_n(f')_{r_n}\to f'$ uniformly on compact subsets of $\Bbb D$ as $n\to\infty$, we have
$$
 \limsup_{n\to\infty} \sup_{\left|\phi(z)\right|\leq r_N}\mu(z)\left|\left(f-f_{r_n}\right)'\left(\phi(z)\right)\right|
 \left|\phi'(z)\right|=0,
$$
where $N\in\Bbb N$ is large enough such that $r_n\geq\frac{1}{2}$ for all $n\geq N$.
Hence we only have to show that
$$
 S:=\limsup_{n\to\infty} \sup_{\left|\phi(z)\right|> r_N}\mu(z)\left|\left(f-f_{r_n}\right)'\left(\phi(z)\right)\right|
 \left|\phi'(z)\right|\leq M_\alpha L.
$$
Indeed, we write $S\leq \limsup_{n\to\infty}\left(S_1+S_2\right)$, where
$$
 S_1=  \sup_{\left|\phi(z)\right|> r_N}\mu(z)\left|f'\left(\phi(z)\right)\right|
 \left|\phi'(z)\right|\hspace{0.2cm}\text{and}\hspace{0.2cm}
 S_2=  \sup_{\left|\phi(z)\right|> r_N}\mu(z)r_n\left|f'\left(r_n\phi(z)\right)\right|
 \left|\phi'(z)\right|.
$$
Then we have
\begin{eqnarray*}
 S_1&=&\sup_{\left|\phi(z)\right|> r_N}\mu(z)\left|f'\left(\phi(z)\right)\right|
 \left|\phi'(z)\right|
 \frac{v_\alpha\left(\phi(z)\right)\left|\sigma'_{\phi(z)}\left(\phi(z)\right)\right|}
 {v_\alpha\left(\phi(z)\right)\left|\sigma'_{\phi(z)}\left(\phi(z)\right)\right|}\\
 &\leq&\frac{4}{\alpha}\|f\|_{\mathcal{B}^\alpha}\sup_{\left|\phi(z)\right|> r_N}\mu(z) \left|\sigma'_{\phi(z)}\left(\phi(z)\right)\right| \left|\phi'(z)\right|\\
 &\leq& \frac{4}{\alpha}\sup_{\left|\phi(z)\right|> r_N}\sup_{|a|>r_N}\mu(z) \left|\sigma'_{a}\left(\phi(z)\right)\right| \left|\phi'(z)\right|
 \leq \frac{4}{\alpha}\sup_{|a|>r_N}\left\|\sigma_a\circ \phi\right\|_{\mathcal{B}^\mu},
\end{eqnarray*}
where we have used the relation (\ref{des-cota-inf-der-sigma}) in the first inequality and the fact that $\left\|f\right\|_{\mathcal{B}^\alpha}\leq 1$ in the second one. Taking limit as $N\to\infty$ we obtain
$$
 \limsup_{n\to\infty}S_1\leq \frac{4}{\alpha} L.
$$
In similar way, we have
\begin{eqnarray*}
 S_2&\leq& \frac{4}{\alpha} \left\|f\right\|_{\mathcal{B}^\alpha} \sup_{\left|\phi(z)\right|> r_N}\mu(z) \left|\sigma'_{\phi(z)}\left(\phi(z)\right)\right| \left|\phi'(z)\right| \frac{r_n v_\alpha\left(\phi(z)\right)}{v_\alpha\left(r_n \phi (z)\right)}\\
 &\leq& \frac{4}{\alpha} \sup_{|a|>r_N}\left\|\sigma_a\circ \phi\right\|_{\mathcal{B}^\mu},
\end{eqnarray*}
since $rv_\alpha(z) < v_\alpha\left(rz\right)$ for all $r\in (0,1)$ and all $z\in\Bbb D$. Therefore
\begin{equation}
  \left\|C_\phi\right\|_e^{\mathcal{B}^\alpha\to\mathcal{B}^\mu}\leq \frac{8}{\alpha}\limsup_{|a|\to 1^-}\left\|\sigma_a\circ \phi\right\|_{\mathcal{B}^\mu}.
\end{equation}
This completes the proof of the theorem.\hfill $\blacksquare$

\vspace{0.2cm}
As an immediate consequence of Theorem \ref{th-principal}, we have the following corollary which generalize a result obtained recently by Malav\'e and Ramos-Fern\'andez in \cite{MR13} and extend a result due to Tjani in \cite{Tj03}. A similar result was found by Gim\'enez, Malav\'e and Ramos-Fern\'andez in \cite{GMR10}, but for the composition operator $C_\phi:\mathcal{B}\to \mathcal{B}^\mu$, where the weight $\mu$  can be extended to non vanishing, complex valued holomorphic function that satisfy a reasonable geometric condition on the Euclidean disk $D(1, 1)$.

\begin{corollary}
 The composition operator $C_\phi$ is compact from $\mathcal{B}^{\alpha}$ into $\mathcal{B}^\mu$  if and only if \, $\phi\in\mathcal{B}^{\mu}$ and
 \begin{equation}\label{limit-compact}
  \lim_{\left|a\right|\to 1^-}\left\|\sigma_a\circ \phi\right\|_{\mu}=0.
 \end{equation}
\end{corollary}


\begin{thebibliography}{00}


\bibitem{At92}  K. Attele.  Toeplitz and Hankel operators on Bergman spaces,  \textit{Hokkaido Math. J.}, {\bf 21}, (1992), 279-293.

\bibitem{CCFR13} R. E. Castillo, D. D. Clahane, J. F. Far\'{\i}as and J. C. Ramos-Fern\'andez. Composition operators from logarithmic Bloch spaces to weighted Bloch spaces, {\it to appear in Appl. Math. Comput.}

\bibitem{ContH00} M. D. Contreras and A. G. Hern\'andez-D\'{\i}az. Weighted composition operators in weighted Banach spaces of analytic functions. \textit{J. Austral Math. Soc. (Serie A)}, {\bf 69}, (2000), 41-60.

\bibitem{CowMac95} C. C. Cowen and B. D. MacCluer. {\em Composition
Operators on Spaces of Analytic Functions}, CRC Press, Boca Raton,
1995.

\bibitem{GMR10} J. Gim\'enez, R. Malav\'e, J. Ramos-Fern\'andez. Composition operators on $\mu$-Bloch type spaces, \textit{Rend. Circ. Mat. Palermo}, {\bf 59}  107-119 (2010).

\bibitem{gblochpaper} O. Hyv\"{a}rinen, M. Kemppainen, M. Lindstr\"{o}m, A. Rautio, and E. Saukko.
 The essential norm of weighted composition operators on weighted Banach spaces of analytic functions, \textit{Integr. Equ. Oper. Theory},
{\bf 72}, (2012), 151-157.


\bibitem{MacZhao03} B. MacCluer and R. Zhao. Essential norms of weighted composition operators between Bloch-type spaces, \textit{Rocky Mountain J. Math.}, {\bf 33},  (2003), 1437-1458.

\bibitem{MM95} K. Madigan, A. Matheson.  Compact composition operators on the Bloch space,  \textit{Trans. Amer. Math. Soc.}, {\bf 347}, (1995), 2679-2687.

\bibitem{MR13} M. Malav\'e Ram\'{\i}rez, J. C. Ramos-Fern\'andez. On a criterion for continuity and compactness of composition operators acting on $\alpha$-Bloch spaces, {\it  C. R. Math. Acad. Sci. Paris} (2013), http://dx.doi.org/10.1016/j.crma.2012.11.013.

\bibitem{MoRo99} A. Montes-Rodr\'{\i}guez. The essential norm of composition operator on Bloch spaces, \textit{Pacific J. Math.}, {\bf 188}, (1999), 339-351.

\bibitem{MoRo00} A. Montes-Rodr\'{\i}guez. Weighted composition operators on weighted Banach spaces of analytic functions, \textit{J. London Math. Soc. (2)}, {\bf 61}, (2000), 872-884.

\bibitem{Ra11EM} J. C. Ramos-Fern\'andez. Composition operators between $\mu$-Bloch spaces. {\it Extracta Math.} {\bf 26}, no. 1, (2011), 75-88.

\bibitem{Tj} M. Tjani.   \textit{Compact composition operators on some M\"{o}bius invariant Banach space}, Ph. D. dissertation, Michigan State University, 1996.

\bibitem{Tj03} M. Tjani. Compact composition operators on Besov spaces. {\it Trans. Amer. Math. Soc.} {\bf 355}, no. 11, 4683-4698, (2003).

\bibitem{Wu07} H. Wulan. Compactness of composition operators on $BMOA$ and $VMOA$. {\it Sci. China Ser. A},  {\bf 50}, no. 7, 997-1004  (2007).

\bibitem{WZheZhu09}  H. Wulan, D.  Zheng, and K. Zhu.  Compact composition operators on $BMO$ and the Bloch space, \textit{Proc. Amer. Math. Soc.}, {\bf 137}, (2009), 3861-3868.

\bibitem{Xi01}  J. Xiao.  Composition operators associated with Bloch-type spaces, \textit{Complex Variables Theory Appl.}, {\bf 46}, (2001), 109-121.

\bibitem{YaZho10} K. Yang, Z. Zhou. Essential norm of the difference of composition operators on Bloch space, \textit{Czechoslovak Math. J.}, {\bf 60(135)},  (2010), 1139-1152.

\bibitem{Zhao10} R. Zhao. Essential norms of composition operators between Bloch type spaces, \textit{Proc. Amer. Math. Soc.}, {\bf 138}, (2010), 2537-2546.
\end{thebibliography}
\end{document}